ArXiv:math/0504058v3 [math.ST] 14 Aug 2007

# MINIMAX AND ADAPTIVE ESTIMATION OF THE WIGNER FUNCTION IN QUANTUM HOMODYNE TOMOGRAPHY WITH NOISY DATA

By Cristina Butucea, Mădălin Guţă[1] and Luis Artiles

*Université Paris X, University of Utrecht and Eurandom*

We estimate the quantum state of a light beam from results of quantum homodyne measurements performed on identically prepared quantum systems. The state is represented through the Wigner function, a generalized probability density on $\mathbb{R}^2$ which may take negative values and must respect intrinsic positivity constraints imposed by quantum physics. The effect of the losses due to detection inefficiencies, which are always present in a real experiment, is the addition to the tomographic data of independent Gaussian noise.

We construct a kernel estimator for the Wigner function, prove that it is minimax efficient for the pointwise risk over a class of infinitely differentiable functions, and implement it for numerical results. We construct adaptive estimators, that is, which do not depend on the smoothness parameters, and prove that in some setups they attain the minimax rates for the corresponding smoothness class.

**1. Introduction.** In 1932 Wigner published a seminal paper [30] in which he introduced a fundamental tool for quantum mechanics known these days as the Wigner function. Glauber extended such techniques to quantum optics where phase space representations of quantum states play an important role in detecting quantum effects in light [7, 13].

*Quantum homodyne tomography* (QHT) is a technique for reconstructing the state of a quantum system from measurement data. It was theoretically proposed in [29] and put in practice for the first time by Smithey et al. [26].

Received April 2005; revised June 2006.

[1]Supported by European IST Programme "Resources for Quantum Information" (RESQ) Contract IST-2001-37559, and by the Netherlands Organisation for Scientific Research (NWO).

This method allows quantum opticians to visualize the Wigner function of newly created states of light and verify whether the theoretical predictions agree with the statistical findings. We mention a few experiments such as the creation of squeezed states [5] and of single-photon-added coherent states [31].

Various aspects of the corresponding ill-posed inverse problem have been analyzed in [9, 23] and [22], and different estimation methods have been proposed by Banaszek et al. [3] and Lvovsky [24]. For an overview of the QHT problem in quantum optics we refer to [21] and for more recent developments to [25].

This paper addresses the statistical problem of estimating the Wigner function of a beam of light from results of QHT measurements on independent, identically prepared beams.

One way to think about quantum tomography as a statistical problem is as follows: the unknown parameter is a joint density $W$ of two variables, $Q$ and $P$. We observe the random variable $(X, \Phi) = (\cos(\Phi)Q + \sin(\Phi)P, \Phi)$ where $\Phi$ is chosen independently of $(Q, P)$, and uniformly in the interval $[0, \pi]$. The joint density of $(X, \Phi)$ can be expressed mathematically in terms of the joint density $W$ of $(Q, P)$, which is allowed to take negative as well as positive values, subject to certain restrictions which guarantee that $(X, \Phi)$ does have a proper probability density. In an ideal situation $W$ would be a density function and then the statistical problem would be to estimate $W$ from independent samples of $(X, \Phi)$. In the context of positron emission tomography this problem has been addressed in [8], which provides minimax rates for the pointwise risk on a class of "very smooth" probability densities. The quantum tomography version where $W$ is a proper Wigner function is treated along similar lines in [16] with the important difference that the proof of the lower bound requires the construction of a "worst parametric family" of Wigner functions rather than probability densities.

In this paper we consider a statistical problem which is more relevant for the experimentalist confronted with various noise sources corrupting the ideal data $(X, \Phi)$. It turns out that a good model for a realistic quantum tomography measurement amounts to replacing $(X, \Phi)$ by the noisy observations $(Y, \Phi)$, where $Y := \sqrt{\eta}X + \sqrt{(1 - \eta)/2}\xi$, with $\xi$ a standard Gaussian random variable independent of $(X, \Phi)$. The parameter $0 < \eta < 1$ is called the detection efficiency and represents the proportion of photons which are not detected due to the losses in the measurement process. This is the statistical problem of this paper, a combination of two classical problems: noise deconvolution and PET tomography. The nonclassical feature is that although all the one-dimensional projections of $W$ are indeed *bona fide* probability densities, the underlying two-dimensional "joint density" need not itself be a *bona fide* joint probability density, but can have small patches of "negative probability."

So far there has been little attention paid to this problem by statisticians, although on the one hand it is an important statistical problem coming up in modern physics, and on the other hand it is "just" a classical nonparametric statistical inverse problem. A first step in the direction of estimating $\rho$ has been made in [2], where consistency results are presented for linear and sieve maximum likelihood estimators. We recommend this paper as a complement to the present one.

Section 2 starts with a short introduction to quantum mechanics followed by the particular problem of estimating the Wigner function in quantum homodyne tomography. In Section 2.3 we describe some features of Wigner functions and show to what extent these functions differ from probability densities on the plane. The section ends with a description of the experimental set-up and the derivation of the Gaussian noise from physical principles.

Section 3 contains the main results of this paper. We assume that the unknown Wigner function belongs to a class $\mathcal{A}(\beta, r, L)$ of "very smooth" functions similar to those of [6, 8] and [16]. The estimator has a standard kernel-type form performing in one step the deconvolution and the inverse Radon transform. In Proposition 1 we compute upper bounds for the pointwise risk. Theorem 1 establishes the lower bound and gives the minimax rate, which is slower than any power of $1/n$ but faster than any power of $1/\log n$. Rates with a similar behavior have been obtained in [6], which inspired some of the results obtained in this paper. Adaptive estimators can be derived in some cases (when $r \leq 1$) (see Theorem 2), converging at the same rates as their nonadaptive correspondents.

In Section 4 we present results of computer simulations for a few quantum states, among which is the Schrödinger cat state which is expected to be produced in the lab in the future. Section 5 collects the proof of Proposition 1 and a sketch of the proof of the adaptive upper bounds.

Section 6 concentrates on the proof of the lower bound for the pointwise risk. For this we construct a pair of Wigner functions $W_{1,2}$ belonging to the class $\mathcal{A}(\beta, r, L)$ such that the distance between them is large enough and the $\chi^2$ distance between the likelihoods of the corresponding models is small. It is now a well-known lower-bounds principle that the best rate of estimation can be viewed as the largest distance between parameters in order to detect the change in the statistical model. This construction is original in the statistics literature as it relies on the positivity of the corresponding density matrices $\rho_1$ and $\rho_2$ rather than of the Wigner functions themselves.

**2. Physical background of quantum tomography.** In this section we present a short introduction to quantum mechanics in as far as it is needed for understanding the background of our statistical problem. The reader who is not interested in the physics can skip this section and continue with Section 3. In Section 2.2 we describe the measurement technique called quantum homodyne tomography and show how this can be used to estimate the Wigner function which is a particular parametrization of the quantum state of a monochromatic pulse of light. More details on Wigner functions can be found in Section 2.3. The main issue tackled in this paper is the influence of noise due to the detection process on the estimation of the Wigner function. The experimental setup of quantum homodyne tomography with noisy observations is discussed in Section 2.4.

For more background material we refer to the textbook [21] on quantum optics and quantum tomography, the paper [2] which deals with the problem of quantum tomography from a statistical perspective, the review paper on quantum statistical inference [4] and the classic textbooks on quantum statistics [17] and [18].

2.1. *Short excursion into quantum mechanics.* Quantum mechanics is the theory which describes the physical phenomena taking place at the microscopic level such as the emission and absorption of light by individual atoms, the detection of light photons. As a theory about physical reality, quantum mechanics makes predictions about the results of measurements performed in the lab. Such predictions are statistical in nature in the sense that in general we cannot infer the result of a measurement on a single quantum system but only the probability distribution of results of identical measurements performed on a statistical ensemble of identically prepared systems. Any such distribution is a function of the state in which the system is prepared, and of the performed measurement. Our statistical problem can then be briefly described as follows: estimate the state based on results of measurements on a number of identically prepared systems.

Mathematically, the main concepts of quantum mechanics are formulated in the language of self-adjoint operators acting on Hilbert spaces. The reader who is not familiar with this theory may think of finite-dimensional Hilbert spaces $\mathbb{C}^d$, and $d \times d$ matrices as operators on $\mathbb{C}^d$. To every quantum system one can associate a complex Hilbert space $\mathcal{H}$ with inner product $\langle \cdot, \cdot \rangle$ whose vectors represent the wave functions of the system or pure states, as we will see below. In general, a state is described by a *density matrix*, which is a compact operator $\rho$ on $\mathcal{H}$ having the following properties:

1. Self-adjoint: $\rho = \rho^*$, where $\rho^*$ is the adjoint of $\rho$.
2. Positive: $\rho \geq 0$, or equivalently $\langle \psi, \rho\psi \rangle \geq 0$ for all $\psi \in \mathcal{H}$.
3. Trace 1: $\mathrm{Tr}(\rho) = 1$.

The positivity property implies that all the eigenvalues of $\rho$ are nonnegative, and by the trace property, they sum up to 1. The reader may have noticed that the above requirements are reminiscent of the properties of probability distributions, and this connection will be strengthened in a moment when we discuss the distribution of measurement results.

Before that we will take a look at the structure of the space of states on a given Hilbert space $\mathcal{H}$. Clearly, the convex combination $\lambda \rho_1 + (1 - \lambda)\rho_2$ of two density matrices $\rho_1$ and $\rho_2$ is a density matrix again and it corresponds to the state obtained as the result of randomly performing one of the two preparation procedures with probabilities $\lambda$ and, respectively, $1 - \lambda$. The extremals of the convex set of states are called *pure states* and are represented by one-dimensional orthogonal projection operators. Indeed an arbitrary density matrix can be brought to the diagonal form

$$\rho = \sum_{i=1}^{\dim \mathcal{H}} \lambda_i \mathbf{P}_i,$$

where $\mathbf{P}_i$ is the projection onto the one-dimensional space generated by the eigenvector $e_i \in \mathcal{H}$ of $\rho$ and $\lambda_i \geq 0$ is the corresponding eigenvalue, that is, $\rho e_i = \lambda_i e_i$.

The predictions made by quantum mechanics can be tested in the lab by performing measurements on quantum systems. We will now give the mathematical description of a measurement with space of outcomes given by the measure space $(\Omega, \Sigma)$. If the system is prepared in the state $\rho$, then the result is *random* and has probability distribution $P_\rho$ over $(\Omega, \Sigma)$ such that the map $\rho \mapsto P_\rho$ is affine, that is, it maps a convex combination of states into the corresponding convex combination of probability distributions. This can be naturally interpreted as saying that for any mixed state $\lambda \rho_1 + (1 - \lambda)\rho_2$, the distribution of the results will reflect the randomized preparation procedure.

The most common measurement is that of an observable such as energy, position, spin, and so on. An observable is described by a *self-adjoint operator* $\mathbf{X} = \mathbf{X}^*$ on the Hilbert space $\mathcal{H}$ and we suppose here for simplicity that it has a discrete spectrum, that is, it can be written in the diagonal form

$$(1) \qquad \mathbf{X} = \sum_{a=1}^{\dim \mathcal{H}} x_a \mathbf{P}_a,$$

with $x_a \in \mathbb{R}$ the eigenvalues of $\mathbf{X}$, and $\mathbf{P}_a$ one-dimensional projections onto the eigenvectors of $\mathbf{X}$. The result of the measurement of the observable $\mathbf{X}$ will be denoted by $X$ and is a random variable with values in the set $\Omega = \{x_1, x_2, \ldots\}$. When the system is prepared in the state $\rho$, the result $X$ has the distribution

$$(2) \qquad \mathbb{P}_\rho[X = x_a] = \mathrm{Tr}(\mathbf{P}_a \rho).$$

Notice that the conditions defining the density matrices insure that $\mathbb{P}_\rho$ is indeed a probability distribution. In particular, the expectation on $X$ in the state $\rho$ is

$$(3) \qquad \mathbb{E}_\rho[X] := \sum_{a=1}^{\dim \mathcal{H}} x_a \mathbb{P}_\rho[X = x_a] = \mathrm{Tr}(\mathbf{X}\rho),$$

and the characteristic function is given by

$$(4) \qquad \mathbb{E}_\rho[\exp(itX)] = \mathrm{Tr}[\exp(it\mathbf{X})\rho].$$

Measurements with continuous outcomes as well as outcomes in an arbitrary measure space can be described in a similar way by using the spectral theory of self-adjoint operators [18].

Suppose that a preparation procedure produces an unknown state $\rho$. It is clear that in general no individual measurement can completely determine the state but only gives us statistical information about $\mathbb{P}_\rho$ and thus indirectly about $\rho$. The problem of state estimation should then be considered in the context of measurements on a large number of systems which are identically prepared in the state $\rho$. Here we consider the simplest situation when we perform identical and independent measurements on each of the $n$ systems separately.

2.2. *Quantum homodyne tomography and the Wigner function.* The statistical problem analyzed in this paper is that of estimating a function $W_\rho \colon \mathbb{R}^2 \to \mathbb{R}$ from i.i.d. data $(Y_1, \Phi_1), \ldots, (Y_n, \Phi_n)$ with distribution $\mathbb{P}_\rho^\eta$ on $\mathbb{R} \times [0, \pi]$. In this subsection we will give an account of the physical origin of this problem.

The quantum system is monochromatic light in a cavity, whose state is described by (infinite-dimensional) density matrices on the Hilbert space of complex-valued square integrable functions on the line $\mathbb{L}_2(\mathbb{R})$. The function of interest $W_\rho$ is called the *Wigner function* and depends in a one-to-one fashion on the state $\rho$ of the light.

Two important observables of this quantum system are the electric and magnetic fields whose corresponding self-adjoint operators on $\mathbb{L}_2(\mathbb{R})$ are given by

$$\mathbf{Q}\psi(x) = x\psi(x) \quad \text{and, respectively,} \quad \mathbf{P}\psi(x) = -i\frac{d\psi}{dx}.$$

The Wigner function $W_\rho \colon \mathbb{R}^2 \to \mathbb{R}$ is much like a joint probability density for these quantities; for instance, its marginals along any direction $\phi \in [0, \pi]$ in the plane which are given by the *Radon transform* of $W_\rho$,

$$(5) \qquad \mathcal{R}[W_\rho](x, \phi) = \int_{-\infty}^{\infty} W_\rho(x\cos\phi - t\sin\phi, x\sin\phi + t\cos\phi)\, dt,$$

are *bona fide* probability densities and correspond to the measurement of the *quadrature* observables $\mathbf{X}_\phi := \mathbf{Q}\cos\phi + \mathbf{P}\sin\phi$. However, in quantum mechanics noncommuting observables such as $\mathbf{Q}$ and $\mathbf{P}$ *cannot* be measured simultaneously; thus we cannot speak of their joint probability distribution. This fact is reflected at the level of the Wigner function, which need not be positive; indeed, it might contain patches of "negative probability."

Thus, for a given quantum system prepared in state $\rho$ we can measure only one of the quadratures $\mathbf{X}_\phi$ for some phase $\phi$ and we obtain a result with probability density $p_\rho(x|\phi) = \mathcal{R}[W_\rho](x, \phi)$. Let us consider now that we have $n$ quantum systems prepared in the same state $\rho$ and we measure the quadrature $\mathbf{X}_{\Phi_i}$ on the $i$th system with phases $\Phi_i$ chosen independently with uniform distribution on $[0, \pi]$. We obtain independent identically distributed results $(X_1, \Phi_1), \dots, (X_n, \Phi_n)$ with density $p_\rho(x, \phi) = p_\rho(x|\phi)$ with respect to the measure $\frac{1}{\pi}\lambda$, where $\lambda$ is the Lebesgue measure on $\mathbb{R} \times [0, \pi]$. The Radon transform $\mathcal{R}: W_\rho \mapsto p_\rho(x, \phi)$ is well known in statistics for its role in tomography problems such as positron emission tomography (PET) [28], and has a broad spectrum of other applications ranging from astronomy to geophysics [10]. In PET one estimates a probability density $f$ on $\mathbb{R}^2$ related to the tissue distribution in a cross section of the human body from i.i.d. observations $(X_1, \Phi_1), \dots, (X_n, \Phi_n)$, with probability density equal to $\mathcal{R}[f]$. The observations are obtained by recording events whereby pairs of photons emitted at the collision of a positron and an electron hit detectors placed in a ring around the body after flying in opposite directions along an axis determined by an angle $\phi \in [0, \pi]$. The difference with our situation is that the role of the unknown distribution is played by the Wigner function, which as we mentioned is not necessarily positive in the usual sense but carries an intrinsic positivity constraint in the sense that it corresponds to a density matrix (see Section 2.3). Another difference with respect to PET is that in QHT the experimenter can decide how to choose the phases $\Phi_i$. Indeed, in some experiments the phases are equidistant, that is, they take one of the values $\frac{l}{k}\pi$ where $l$ runs from 0 to $k - 1$ for some $k \in \mathbb{N}$, but one has now the additional problem of how to choose $k$ as a function of $n$. We believe that by using uniformly distributed phases one does not incur any loss in the asymptotic rates, but it remains an interesting open question whether a specially designed choice of phases can improve the results. This may be the case for some parametric classes of Wigner functions with an asymmetric aspect like those corresponding to squeezed states (see Section 2.3).

2.3. *Properties of Wigner functions.* The physics literature on Wigner functions and other types of "phase space functions" is vast, but a starting point for the interested reader may be the book [21]. Here we focus on the similarities and the differences with usual probability densities encountered in PET.

Consider the space of Hilbert–Schmidt operators on $\mathbb{L}_2(\mathbb{R})$,

$$\mathcal{T}_2 := \{\mathbf{A} \in \mathcal{B}(\mathbb{L}_2(\mathbb{R})) : \|\mathbf{A}\|_2^2 = \mathrm{Tr}(\mathbf{A}^*\mathbf{A}) < \infty\},$$

on which there exists an inner product $\langle \mathbf{A}, \mathbf{B} \rangle_2 = \mathrm{Tr}(\mathbf{A}^*\mathbf{B})$, and notice that the density matrices form a closed subset of $\mathcal{T}_2$. The Wigner function $W_\mathbf{A}$ is the image of **A** through the linear map $W : \mathcal{T}_2 \to \mathbb{L}_2(\mathbb{R}^2)$ defined by the property that the Fourier transform $\mathcal{F}_2$ with respect to *both* variables has the expression

$$(6) \qquad \widetilde{W}_{\mathbf{A}}(u, v) := \mathcal{F}_2[W_{\mathbf{A}}](u, v) = \mathrm{Tr}(\mathbf{A} \exp(iu\mathbf{Q} + iv\mathbf{P})).$$

In particular, this defines the Wigner function $W_\rho$ of the state with density matrix $\rho$. By passing to the polar coordinates $(u, v) = (t \cos\phi, t \sin\phi)$ we have $u\mathbf{Q} + v\mathbf{P} = t\mathbf{X}_\phi$, and using (4) together with the fact that $p_\rho(\cdot | \phi)$ is the density for measuring $\mathbf{X}_\phi$ we have

$$(7) \qquad \widetilde{W}_\rho(u, v) = \mathrm{Tr}(\rho \exp(it\mathbf{X}_\phi)) = \mathcal{F}_1[p_\rho(\cdot | \phi)](t),$$

where the Fourier transform $\mathcal{F}_1$ in the last term is with respect to the first variable, keeping $\phi$ fixed. The reader familiar with PET may recognize that the composition $\mathcal{F}_2 \circ \mathcal{F}_1$ mapping $p_\rho$ into $W_\rho$ is just the inverse Radon transform [10], proving our assertion that QHT is about the tomography of the Wigner function.

It can be shown that the map $W : \mathcal{T}_2 \to \mathbb{L}_2(\mathbb{R}^2)$ is isometric up to a constant:

$$(8) \qquad \langle \mathbf{A}, \mathbf{B} \rangle_2 = 2\pi \langle W_{\mathbf{A}}, W_{\mathbf{B}} \rangle := 2\pi \iint \overline{W}_{\mathbf{A}}(q, p) W_{\mathbf{B}}(q, p) \, dq \, dp,$$

and this fact is often used as a tool for calculating the expectation of an observable $\mathbf{X} \in \mathcal{T}_2$ similarly to the way it is done in classical probability:

$$(9) \qquad \mathrm{Tr}(\rho \mathbf{X}) = 2\pi \iint W_{\mathbf{X}}(q, p) W_\rho(q, p) \, dq \, dp.$$

Let us come back to our physical system, the light in a cavity, and consider its energy, which is given by the sum of intensities of the electric and magnetic fields $\mathbf{H} := \frac{1}{2}(\mathbf{Q}^2 + \mathbf{P}^2)$. As predicted by Einstein before the creation of quantum theory, the possible values that this observable may take are "quantized," which can be explained if we think of light as a packet of photons with each photon contributing a fixed quantum of energy. Indeed, by solving the eigenvalue problem we find $\mathbf{H}\psi_j = (j + 1/2)\psi_j$ where $\{\psi_j\}_{j \geq 0}$ is an orthonormal basis of $\mathbb{L}_2(\mathbb{R})$ whose vectors have the physical interpretation of pure states with precisely $j$ photons and are given by

$$(10) \qquad \psi_j(x) = \frac{1}{\sqrt{\sqrt{\pi} 2^j j!}} H_j(x) e^{-x^2/2},$$

where $H_j(x)$ are the Hermite polynomials (see, e.g., [12]).

Notably, the vacuum state corresponding to zero photons has nonzero energy $1/2$, a purely quantum phenomenon called *vacuum fluctuations* reflected in the fact that the distributions of $\mathbf{Q}$ and $\mathbf{P}$ are Gaussian with variance $1/2$. We would like to stress here that the Gaussian distribution emerges directly from physical principles and it is the same Gaussian character of the vacuum which will lead to our model for the detection noise in Section 2.4.

An interesting consequence of relation (9) is found by taking $\mathbf{X}$ to be the vacuum state $\mathbf{P}_{\psi_0}$ whose Wigner function is $W_{\mathbf{X}}(q,p) = \exp(-q^2 - p^2)/\pi$. Then, as the left-hand side of the equation is positive, this implies that the negative patches of $W_\rho$ around the origin must be balanced by positive ones in such a way that the integral remains positive. In fact this property holds for any point in the plane and the localized oscillations of the Wigner function are a signature of nonclassical states, such as states with a fixed number of photons or the so-called "Schrödinger cat states" like the one estimated in Figure 3.

On the other hand, there exist probability densities that are not Wigner functions, for example, the latter cannot be too "peaked" (cf. [21]):

$$(11) \qquad |W_\rho(q,p)| \le \frac{1}{\pi} \qquad \text{for all } (q,p) \in \mathbb{R}^2.$$

A general density matrix $\rho$ can be seen as an infinite-dimensional matrix with coefficients $\rho_{jk} = \langle \psi_j, \rho \psi_k \rangle$ for $j, k \ge 0$ such that $\sum_{k \ge 0} \rho_{kk} = 1$ (trace 1), and $[\rho_{jk}] \ge 0$ (positive definite matrix). In particular, the diagonal elements $p_k = \rho_{kk}$ represent the probability of measuring $k$ photons for a system in state $\rho$. The density $p_\rho(x, \phi)$ is given in terms of the matrix elements of $\rho$ by

$$(12) \quad p_\rho(x, \phi) = \frac{1}{\pi} \sum_{j,k=0}^{\infty} \rho_{jk} p_{jk}(x, \phi) := \frac{1}{\pi} \sum_{j,k=0}^{\infty} \rho_{jk} \psi_j(x) \psi_k(x) e^{-i(j-k)\phi},$$

and a similar formula holds for the Wigner function $W_\rho(q,p) = \sum_{j,k=0}^{\infty} \rho_{jk} W_{jk}(q, p)$, with $W_{jk}$ such that $\mathcal{R}[W_{jk}] = p_{jk}$. For any density matrices $\rho, \tau$ (8) can be written

$$\|W_\rho - W_\tau\|_2^2 := \iint |W_\rho(q,p) - W_\tau(q,p)|^2 \, dp \, dq$$

$$(13)$$

$$= \frac{1}{2\pi} \|\rho - \tau\|_2^2 := \frac{1}{2\pi} \sum_{j,k=0}^{\infty} |\rho_{jk} - \tau_{jk}|^2.$$

Some examples of quantum states that can be created at this moment in the lab are given in Table 1 of [2]. Typically, the corresponding Wigner functions have a Gaussian tail but need not be positive. As a consequence of (11) not all two-dimensional Gaussian distributions are Wigner functions, but only those for which the determinant of the covariance matrix is at least $\frac{1}{4}$. Equality is obtained for a remarkable set of states called *squeezed vacuum* states having Wigner functions $W(q, p) = \frac{1}{\pi} \exp(-e^{2\xi} q^2 - e^{-2\xi} p^2)$, determined by the squeezing factor $\xi$. More generally, the celebrated *Heisenberg uncertainty relation* says that for any state $\rho$ the noncommuting observables $\mathbf{P}$ and $\mathbf{Q}$ cannot have probability distributions such that the product of their variances is smaller than $\frac{1}{4}$.

2.4. *Experimental setup and noisy observations.* The optical setup sketched in Figure 1 consists of an additional laser of high intensity $|z|^2 \gg 1$ called a local oscillator, a beam splitter through which the cavity pulse prepared in state $\rho$ is mixed with the laser, and two photodetectors each measuring one of the two beams and producing currents $I_{1,2}$ proportional to the number of photons. An electronic device produces the result of the measurement by taking the difference of the two currents, integrating it over the time interval of the pulse, and rescaling it by a factor proportional to $|z|$ (see below). A detailed analysis taking into account various losses (mode mismatching, detection inefficiency) in the detection process can be found in [21]. It turns out that all these losses can be modeled by a Gaussian noise in the measurement results, and here we detail only the case of detection inefficiency. In the high photon number regime $|z|^2 \gg 1$ the (integrated) current depends linearly on the intensity of the beam with a proportion $\eta < 1$ of the photons being detected. The process can be described classically by considering that each individual photon has probability $\eta$ of being detected and $1-\eta$ of being absorbed without detection. Thus in a beam of $j$ photons the probability of detecting $k \leq j$ is $b_k^j(\eta) = \binom{j}{k}\eta^k(1-\eta)^{j-k}$, and for an incoming state $\rho$ we obtain the probability distribution of the results $p_k(\eta) = \sum_{j=k}^{\infty} \rho_{jj} b_k^j(\eta)$. This "photon lottery" can be equivalently described by replacing the realistic detector with an ideal one in front of which we place an imaginary beam splitter (see Figure 1) which has transmissivity $t = \sqrt{\eta}$ and reflectivity $r = \sqrt{1-\eta}$.

In order to understand why this is the case and how the measurement noise appears, we will present two equivalent pictures of the action of the beam splitter stemming from the *wave-particle duality* typical in quantum mechanics. As shown in Figure 1 a beam splitter receives two incoming beams and has two outgoing beams as output. In the case of the imaginary beam splitter sitting in front of the detector, one of the incoming beams is the vacuum and let us assume that the beam to be measured has $j$ photons. Then the joint state of the two beams is $\psi_0 \otimes \psi_j \in \mathbb{L}_2(\mathbb{R}) \otimes \mathbb{L}_2(\mathbb{R})$ and the transformation to the outgoing vector is $\psi_0 \otimes \psi_j \mapsto \sum_{k=0}^{j} [b_k^j(\eta)]^{1/2} \psi_{j-k} \otimes \psi_k$, which simply means that with probability $b_k^j(\eta)$ we get $k$ photons going to the ideal detector and $j-k$ will not be detected, as described above.

The second description is in terms of the transformation of the electric and magnetic field operators of the beams denoted by $(\mathbf{Q}_l, \mathbf{P}_l)$ and $(\mathbf{Q}_r, \mathbf{P}_r)$, with the first couple acting on the left side of the tensor product $\mathbb{L}_2(\mathbb{R}) \otimes \mathbb{L}_2(\mathbb{R})$ and the second pair on the right side. The fields of the outgoing beams are $\mathbf{Q}'_l = t\mathbf{Q}_l - r\mathbf{Q}_r$, $\mathbf{Q}'_r = r\mathbf{Q}_l + t\mathbf{Q}_r$ and similarly for $\mathbf{P}$'s.

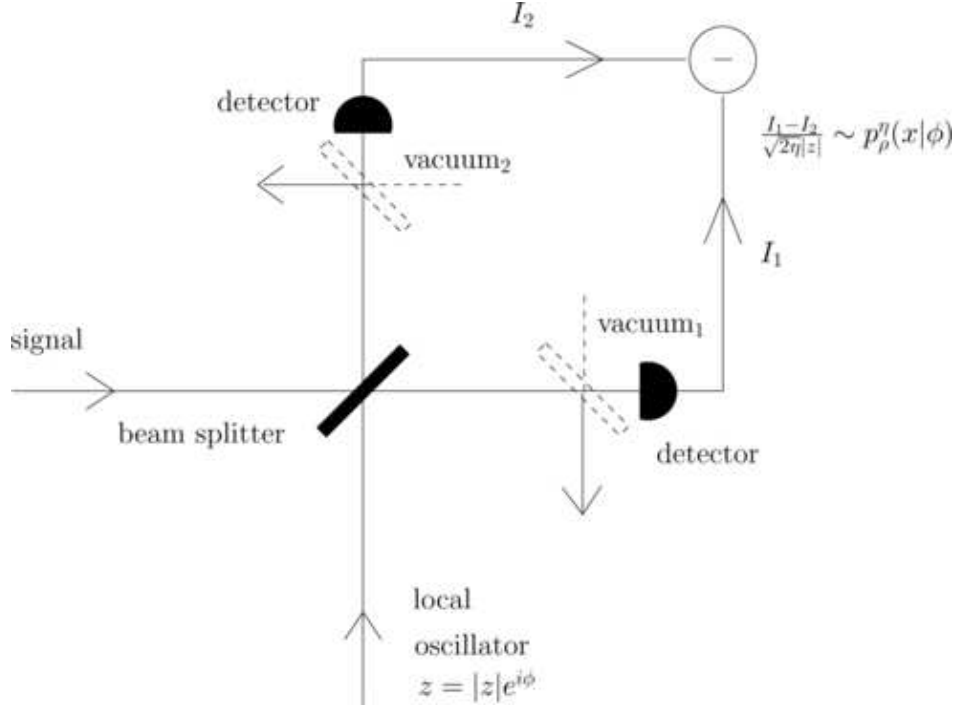

Fig. 1. *Quantum homodyne tomography measurement setup.*

Then by computing the combined effects of the beam splitters, we have the fields arriving at the two detectors, $\mathbf{Q}_1 = \frac{t}{\sqrt{2}}[\mathbf{Q} + \mathbf{Q}_{\mathrm{lo}}] - r\mathbf{Q}_{1\mathrm{vac}}$ and $\mathbf{Q}_2 = \frac{t}{\sqrt{2}}[\mathbf{Q} - \mathbf{Q}_{\mathrm{lo}}] - r\mathbf{Q}_{2\mathrm{vac}}$, and similarly for $\mathbf{P}_1, \mathbf{P}_2$. We remind the reader that the number of photons in a beam is described by $\mathbf{N} := \frac{1}{2}(\mathbf{Q}^2 + \mathbf{P}^2 - \mathbf{1})$. Using the fact that in the limit $|z|^2 \gg 1$ the laser can be treated classically by replacing $\mathbf{Q}_{\mathrm{lo}}$ by $\frac{|z|}{\sqrt{2}}\cos\phi$ and $\mathbf{P}_{\mathrm{lo}}$ by $\frac{|z|}{\sqrt{2}}\sin\phi$, we get

$$\mathbf{N}_1 - \mathbf{N}_2 = \sqrt{2}t|z|[(t\mathbf{Q}_\phi + r\mathbf{Q}_\phi^{\mathrm{vac}}) + O(|z|^{-1})],$$

with $O(|z|^{-1})$ a term whose variance is bounded by $C/|z|$, and $\mathbf{Q}_\phi^{\mathrm{vac}}$ a quadrature operator of a vacuum mode accounting for the two fictitious beam splitters. Thus in the limit $|z| \to \infty$ the rescaled integrated current difference $I_1 - I_2/\sqrt{2\eta}|z|$ has the same distribution as $t\mathbf{Q}_\phi + r\mathbf{Q}_\phi^{\mathrm{vac}}$, that is, that of the sum of two independent random variables $Y := \sqrt{\eta}X + \sqrt{(1-\eta)/2}\xi$, where $X \sim P_\rho(\cdot|\phi)$ is the result of measuring $X_\phi$, $\xi$ has the $N(0,1)$ law and $\frac{1}{\sqrt{2}}\xi$ has the distribution of the quadrature in the vacuum (see Section 2.3).

The efficiency-corrected probability density is then the convolution

$$
(14) \quad p_\rho^\eta(y, \phi) = (\pi(1-\eta))^{-1/2} \int_{-\infty}^{\infty} p_\rho(x, \phi) \exp\left[-\frac{\eta}{1-\eta}(x - \eta^{-1/2}y)^2\right] dx.
$$

Finally, the constants $|z|$ and $\eta$ are measured in advance as part of the calibration of the experiment and are considered to be known.

## 3. Statistical procedure and results.
For convenience we summarize now the statistical problem tackled in this paper.

Consider $(X_1, \Phi_1), \ldots, (X_n, \Phi_n)$, independent identically distributed random variables with values in $\mathbb{R} \times [0, \pi]$ and distribution $\mathbb{P}_\rho$ having density $p_\rho(x, \phi)$ with respect to $\frac{1}{\pi}\lambda$, $\lambda$ being the Lebesgue measure on $\mathbb{R} \times [0, \pi]$, given by

$$
p_\rho(x, \phi) = \mathcal{R}[W_\rho](x, \phi),
$$

where $\mathcal{R}$ is the Radon transform defined in (5) and $W_\rho : \mathbb{R}^2 \to \mathbb{R}$ is a so-called Wigner function which we want to estimate. The space of all possible Wigner functions is parametrized by infinite-dimensional matrices $\rho = [\rho_{jk}]_{j,k=0}^{\infty}$ such that $\operatorname{Tr}\rho = 1$ (trace 1) and $\rho \geq 0$ (positive definite), in the way indicated by (6). Moreover, the correspondence between $\rho$ and $W_\rho$ is one-to-one and isometric with respect to the $\mathbb{L}_2$ norms as in (13). The properties of Wigner functions have been discussed in Section 2.3, in particular the fact that $W_\rho$ may take negative values.

What we observe are not the variables $(X_\ell, \Phi_\ell)$ but the noisy ones $(Y_1, \Phi_1), \ldots, (Y_n, \Phi_n)$, where

$$
(15) \quad Y_\ell := \sqrt{\eta}X_\ell + \sqrt{(1-\eta)/2}\,\xi_\ell,
$$

with $\xi_\ell$ a sequence of independent identically distributed standard Gaussians which are independent of all $(X_j, \Phi_j)$. The parameter $0 < \eta < 1$ is known and we denote by $p_\rho^\eta$ the density of $(Y_\ell, \Phi_\ell)$ given by the convolution (14). The aim is to recover the Wigner function $W_\rho$ from the noisy observations.

*Class of Wigner functions.* In order to apply the minimax estimation technology we will assume that the unknown Wigner function is infinitely differentiable and belongs to the following class described via its Fourier transform:

$$
\mathcal{A}(\beta, r, L) = \left\{ W_\rho \text{ Wigner function} : \int |\widetilde{W}_\rho(w)|^2 e^{2\beta\|w\|^r}\, dw \leq (2\pi)^2 L \right\},
$$

where $0 < r \leq 2$, and $\beta, L > 0$. From now on we denote by $\langle \cdot, \cdot \rangle$ and $\|\cdot\|$ the usual Euclidean scalar product and norm, while $C(\cdot)$ will denote positive constants depending on parameters given in the parentheses. From the physical point of view the choice of a class of very smooth Wigner functions seems to be quite reasonable considering that to date no quantum state of light has been constructed which does not satisfy such conditions. The reason for the difficulty in engineering states with less smooth Wigner functions is that the interactions needed to produce such states should be very nonlinear in the electric and magnetic fields while it is known that photons are rather weakly interacting particles. For example, until recently the creation of squeezed states requiring a quadratic interaction was a not-trivial achievement [5]. We mention here without proof the result of a computation showing that if a density matrix $\rho$ satisfies the condition $\mathrm{Tr}(\rho \exp[a\mathbf{N}^{r/2}]) < \infty$ for some $a, r > 0$, then $W_\rho \in \mathcal{A}(\beta, r, L)$ for some $\beta, L > 0$. In light of the previous argument we consider that this condition is actually rather weak.

*Estimation method.* For the problem of estimating a probability density $f : \mathbb{R}^2 \to \mathbb{R}$ directly from data $(X_\ell, \Phi_\ell)$ with density $\mathcal{R}[f]$ we refer to the literature on X-ray tomography and PET, studied in [8, 19, 20, 28], and the references therein. In the context of tomography of bounded objects with noisy observations, Goldenshluger and Spokoiny [14] solved the problem of estimating the borders of the object (the support). For the problem of Wigner function estimation when no noise is present, we mention the parallel work [16].

Let $N^\eta$ denote the density of the rescaled noise $\sqrt{(1-\eta)/2}\,\xi$ and let $\widetilde{N}^\eta$ be its Fourier transform. Denote by $p_\rho^\eta(y, \phi)$ the probability density of $(Y_\ell, \Phi_\ell)$ in (14). Then

$$p_\rho^\eta(y, \phi) = \int_{-\infty}^{\infty} \frac{1}{\sqrt{\eta}} p_\rho\left(\frac{y-x}{\sqrt{\eta}}, \phi\right) N^\eta(x)\, dx := \left(\frac{1}{\sqrt{\eta}} p_\rho\left(\frac{\cdot}{\sqrt{\eta}}, \phi\right) * N^\eta\right)(y),$$

where $p * q(y) = \int p(y-x)q(x)\, dx$ denotes the convolution of $p$ and $q$. Via a change of variable we can write $p_\rho^\eta(y, \phi)$ as in (14). In the Fourier domain this relation becomes $\mathcal{F}_1[p_\rho^\eta(\cdot, \phi)](t) = \mathcal{F}_1[p_\rho(\cdot, \phi)](t\sqrt{\eta})\widetilde{N}^\eta(t)$, where $\mathcal{F}_1$ denotes the Fourier transform with respect to the first variable.

In this paper we modify the usual tomography kernel in order to take into account the additive noise on the observations and construct a kernel $K_h^\eta$ that asymptotically performs both deconvolution and inverse Radon transformation on our data. Let us define the estimator

(16) $$\widehat{W}_{h,n}^\eta(q, p) = \frac{1}{n} \sum_{\ell=1}^n K_h^\eta\left(q\cos\Phi_\ell + p\sin\Phi_\ell - \frac{Y_\ell}{\sqrt{\eta}}\right),$$

where $0 < \eta < 1$ is a fixed parameter, and the kernel is defined by

(17)
$$K_h^\eta(u) = \frac{1}{4\pi} \int_{-1/h}^{1/h} \frac{\exp(-iut)|t|}{\widetilde{N}^\eta(t/\sqrt{\eta})}\, dt,$$

$$\widetilde{K}_h^\eta(t) = \frac{1}{2} \frac{|t|}{\widetilde{N}^\eta(t/\sqrt{\eta})} I(|t| \le 1/h),$$

and $h > 0$ tends to 0 when $n \to \infty$ in a proper way to be chosen later. For simplicity, let us denote $z = (q, p)$ and $[z, \phi] = q \cos \phi + p \sin \phi$; then the estimator can be written

$$\widehat{W}_{h,n}^{\eta}(z) = \frac{1}{n} \sum_{\ell=1}^{n} K_h^{\eta}\left([z, \Phi_\ell] - \frac{Y_\ell}{\sqrt{\eta}}\right).$$

This is a one-step procedure for treating two successive inverse problems. The main difference with the no-noise problem treated by Guţă and Artiles [16] is that the deconvolution is more difficult than inverse Radon transformation, and thus the techniques for proving the optimality of the method (lower bound) are essentially different. Technically, the no-noise kernel-type estimator has dominating variance, while in the case of noisy observations the bias dominates the variance, as we will see later on.

In Section 3.1 we analyze the mean squared error (MSE) at some fixed point. Our results concern minimax efficiency and adaptive optimality for this problem. We compute an upper bound for the convergence rate of the proposed estimator by minimizing the sum of upper bounds (uniform over the whole class) of the bias and of the variance. The optimality in rate of our estimator follows from the lower bounds, which are proved in Section 6. The meaning of the lower bounds results is that asymptotically, no other estimation technique could outperform our method uniformly over all Wigner functions in the given class. Moreover, we prove the lower bounds, including the asymptotic constant (sharp minimax).

We use a technique based on two hypotheses that appeared in [11] for periodic Sobolev classes and in [6] for classes of supersmooth functions, to which we refer for the details of some of the computations. We concentrate on the main construction involved in the lower bound, that is, the choice of two hypotheses belonging to the fixed class of Wigner functions such that their values in a fixed point are sufficiently different while their corresponding models have likelihoods close to each other.

Despite the generality of a minimax sharp estimator, for practical purposes it is not obvious how to choose the smoothness parameters $r$ and $\beta$. Therefore, an adaptive method (i.e., free of prior knowledge of parameters $\beta$, $r$ and $L$ provided that they are in some set) is designed for classes with $r \leq 1$ in Section 3.2. They behave as well as the previous estimators, provided that we know maximal values of parameters. In particular, this estimator is optimal adaptive (i.e., adaptive and attaining the minimax rate) and efficient. We note that in general such procedures do not always exist. We are fortunate in our case and this is mainly due to the dominating bias.

3.1. *Pointwise estimation.* In this section we give minimax and adaptive results for the pointwise risk (MSE) for the estimator $\widehat{W}_{h,n}^{\eta}$ in (16). The next proposition contains upper bounds for the two components of the risk, the bias and variance, as functions of the parameter $h$ and the number $n$ of samples. The bounds are uniform over all Wigner functions in the class $\mathcal{A}(\beta, r, L)$.

PROPOSITION 1. *Let* $(Y_\ell, \Phi_\ell), \ell = 1, \ldots, n$, *be i.i.d. data coming from the model* (15) *and let* $\widehat{W}_{h,n}^{\eta}$ *be an estimator (with* $h \to 0$ *as* $n \to \infty$*) of the underlying Wigner function* $W_\rho$ *belonging to the class* $\mathcal{A}(\beta, r, L)$, *with* $0 < r \leq 2$. *Then*

$$\sup_{z \in \mathbb{R}^2} \sup_{W_\rho \in \mathcal{A}(\beta, r, L)} |\mathbb{E}[\widehat{W}_{h,n}^{\eta}(z)] - W_\rho(z)|^2 = \frac{L h^{r-2}}{4\pi\beta r} \exp\left(-\frac{2\beta}{h^r}\right)(1 + o(1)),$$

$$\sup_{z \in \mathbb{R}^2} \sup_{W_\rho \in \mathcal{A}(\beta, r, L)} \mathbb{E}[|\widehat{W}_{h,n}^{\eta}(z) - E[\widehat{W}_{h,n}^{\eta}(z)]|^2] \leq \frac{1}{8\gamma^2 n} \exp\left(\frac{2\gamma}{h^2}\right)(1 + o(1)),$$

*where* $\gamma = (1 - \eta)/(4\eta)$, *and* $o(1) \to 0$ *as* $h \to 0$ *and* $n \to \infty$.

The pointwise convergence rate of $\widehat{W}_{h,n}^{\eta}$ with $h = h_{\text{opt}}$ is then shown to be minimax by proving an additional lower bound.

THEOREM 1. *Let* $\beta > 0$, $L > 0$, $0 < r \leq 2$ *and* $(Y_\ell, \Phi_\ell), \ell = 1, \ldots, n$, *be i.i.d. data coming from the model* (15), *and let* $\widehat{W}_{h,n}^{\eta}$ *be as defined in* (16) *with the kernel* $K_h^{\eta}$ *of* (17) *and let the bandwidth* $h_{\text{opt}}$ *be given by the solution of*

$$(18) \qquad \frac{2\beta}{h_{\text{opt}}^r} + \frac{2\gamma}{h_{\text{opt}}^2} = \log n.$$

*Then* $\widehat{W}_{h,n}^{\eta}$ *satisfies the following upper bounds in pointwise distance:*

$$\limsup_{n \to \infty} \sup_{z \in \mathbb{R}^2} \sup_{W_\rho \in \mathcal{A}(\beta, r, L)} \mathbb{E}[|\widehat{W}_{h,n}^{\eta}(z) - W_\rho(z)|^2] \varphi_n^{-2} \leq C,$$

*where the constant* $C$ *and the pointwise rate are*

$$C = 1, \qquad \varphi_n^2 = \frac{L h_{\text{opt}}^{r-2}}{4\pi\beta r} \exp\left(-\frac{2\beta}{h_{\text{opt}}^r}\right) \qquad \text{if } 0 < r < 2,$$

$$C > 0, \qquad \varphi_n^2 = n^{-\beta/(\beta+\gamma)} \qquad\qquad \text{if } r = 2.$$

*Moreover, the previous rate is minimax efficient for* $0 < r < 2$ *and nearly minimax for* $r = 2$; *that is, the following lower bounds hold:*

$$\liminf_{n \to \infty} \inf_{\widehat{W}_n} \sup_{W_\rho \in \mathcal{A}(\beta, r, L)} \mathbb{E}[|\widehat{W}_n(z) - W_\rho(z)|^2] \varphi_n^{-2} \geq 1 \qquad \forall z \in \mathbb{R} \text{ if } 0 < r < 2,$$

$$\liminf_{n \to \infty} \inf_{\widehat{W}_n} \sup_{W_\rho \in \mathcal{A}(\beta, 2, L)} \mathbb{E}[|\widehat{W}_n(z) - W_\rho(z)|^2](n \log n)^{\beta/(\beta+\gamma)} \geq c > 0$$

$$\forall z \in \mathbb{R} \; if \; r = 2,$$

where $\inf_{\widehat{W}_n}$ is taken over all possible estimators $\widehat{W}_n$ of the Wigner function $W_\rho$.

PROOF.    The proof of the lower bounds is given in Section 6.

*Sketch of proof of the upper bounds.* By Proposition 1 we write

$$\sup_{z \in \mathbb{R}^2} \sup_{W_\rho \in \mathcal{A}(\beta, r, L)} \mathbb{E}[|\widehat{W}_{h,n}^\eta(z) - W_\rho(z)|^2] \leq C_B h^{r-2} \exp\left(-\frac{2\beta}{h^r}\right) + \frac{C_V}{n} \exp\left(\frac{2\gamma}{h^2}\right),$$

where $C_B$ and $C_V$ denote the constant terms, depending on $\beta, r, L$ and $\eta$. We select the best bandwidth as $h_{\mathrm{opt}} = \arg\inf_{h>0}\{C_B h^{r-2} \exp(-2\beta/h^r) + C_V/n \exp(2\gamma/h^2)\}$. By taking derivatives we get

$$\frac{2\beta}{h^r} + \frac{2\gamma}{h^2} = \log n + C(1 + o(1)) \qquad \text{as } n \to \infty,$$

where $C > 0$ depends on $\beta, r, L$ and $\eta$. This allows us to take $h_{\mathrm{opt}}$ as in (18) and check that up to constants

$$h_{\mathrm{opt}}^{r-2} \exp\left(-\frac{2\beta}{h_{\mathrm{opt}}^r}\right) = h_{\mathrm{opt}}^{r-2} \cdot \frac{1}{n} \exp\left(\frac{2\gamma}{h_{\mathrm{opt}}^2}\right)(1 + o(1)) \sim h_{\mathrm{opt}}^{r-2} \operatorname{Var}(\widehat{W}_{h_{\mathrm{opt}},n}^\eta(z)),$$

that is, the bias term is asymptotically larger than the variance term, for all $0 < r < 2$, and they are of the same order if $r = 2$.    □

REMARKS ON BANDWIDTHS AND RATES.    The bandwidth (18) and consequently the rates are given in an implicit form. We show now that more explicit expressions can be obtained, if we restrict to values of $r$ in certain intervals.

If $r \leq 1$, then it suffices to take bandwidth

$$h_1 = \left(\frac{\log n}{2\gamma} - \frac{\beta}{\gamma}\left(\frac{\log n}{2\gamma}\right)^{r/2}\right)^{-1/2}$$

and the bias term is larger than the variance term (for $h = h_1$) and of the same order as $\varphi_n^2$ (for $h = h_{\mathrm{opt}}$):

$$\frac{L}{4\pi\beta r}\left(\frac{\log n}{2\gamma}\right)^{1-r/2} \exp\left(-2\beta\left(\frac{\log n}{2\gamma}\right)^{r/2} + o(1)\right).$$

If $1 < r \leq 4/3$, then we take $h_2 = (\frac{\log n}{2\gamma} - \frac{\beta}{\gamma}h_1^{-r})^{-1/2}$ and we get the risk bound (for $h = h_2$) of the same order as $\varphi_n^2$ (for $h = h_{\mathrm{opt}}$):

$$\frac{L}{4\pi\beta r}\left(\frac{\log n}{2\gamma}\right)^{1-r/2} \exp\left(-2\beta\left(\frac{\log n}{2\gamma}\right)^{r/2} + C_1(r, \beta, \gamma)\left(\frac{\log n}{2\gamma}\right)^{r-1} - o(1)\right).$$

In general one has to consider separately the cases $(k-1)/k < r/2 \leq k/(k+1)$.

We deal with a composition of two ill-posed inverse problems with the deconvolution being the dominating factor and the inverse Radon transformation bringing corrections to the usual rates. For $r = 1$ we can compare our result with that of [16] for the idealized tomography model without noise. While the latter is almost parametric, in the presence of deconvolution the rates decrease to a factor $\sqrt{\log n}\exp(-c\sqrt{\log n})$, which is faster than $(\log n)^{-a}$ but slower than power $n^{-a}$ rates, for any $a > 0$. Compared with the density estimation in the convolution model of [6], we get an additional logarithmic factor $h_{\mathrm{opt}}^{-1}/2$ in the rates due to the presence of the inverse Radon transformation. However, as we will see later, an important difference with [6] is the proof of the lower bound requiring the construction of a "most difficult" family of Wigner functions.

3.2. *Optimal adaptive estimation.* In the previous theorem the kernel estimator $\widehat{W}_{h,n}^{\eta}$ has a bandwidth $h = h_{\mathrm{opt}}$ which is the solution of (18) depending on the parameters $\beta$ and $r$ of the class. In the next theorem we will show that there exists an adaptive estimator, that is, not depending on the parameters, performing as well as the former estimators, provided that they lie in the set $\mathcal{B} = \{(\beta, r, L) : \beta > 0, 0 < r < 1, L > 0\}$.

THEOREM 2. *Let* $(Y_\ell, \Phi_\ell), \ell = 1, \ldots, n$, *be i.i.d. data coming from the model* (15). *Then* $\widehat{W}_{h,n}^{\eta}$ *with* $h = h_{\mathrm{ad}}$,

$$h_{\mathrm{ad}} = \left(\frac{2\eta\log n}{1-\eta} - \sqrt{\frac{2\eta\log n}{1-\eta}}\right)^{-1/2},$$

*is an optimal adaptive estimator over the set of parameters* $\mathcal{B}$. *That is, the estimator attains the same upper bounds, for all* $(\beta, r, L) \in \mathcal{B}$,

$$\limsup_{n\to\infty} \sup_{W_\rho \in \mathcal{A}(\beta, r, L)} \mathbb{E}[|\widehat{W}_{h_{\mathrm{ad}},n}^{\eta}(z) - W_\rho(z)|^2]\varphi_n^{-2} \leq 1 \qquad \forall z \in \mathbb{R}^2,$$

*where the rate* $\varphi_n^{-2}$ *is given in Theorem* 1 *for the case* $0 < r < 1$.

For the proof of this theorem we refer to a similar result of [6]. Note that $h_{\mathrm{ad}}$ is a fixed quantity and does not depend on the data. An important consequence is that in conjunction with the lower bounds in Theorem 1, the estimator $\widehat{W}_{h_{\mathrm{ad}},n}^{\eta}$ is optimal adaptive and efficient over the set $\mathcal{B}$ for the pointwise risk. This means it attains the minimax rate and the constant $C = 1$ for an estimator free of $\beta, r$ and $L$ provided that these parameters are in the class $\mathcal{B}$.

**4. Practical implementation.** We study three Wigner functions, each one belonging to some class $\mathcal{A}(\beta, 2, L)$ with arbitrary $\beta < 1/4$. The one- and two-photon states are described by diagonal density matrices with $\rho_{jj} = \delta_{j,1}$ and, respectively, $\rho_{jj} = \delta_{j,2}$, and can be readily produced in the lab. The third state is a so-called Schrödinger cat state which is represented by the sum of two vectors corresponding to laser states, and which may be available experimentally in the near future.

For the one-photon state, we simulated $n = 5000$ noisy data $(Y_\ell, \Phi_\ell)$ by first generating $(X_\ell, \Phi_\ell)$ having density $p_\rho(x, \phi)$ and then adding the noise by using standard Gaussians $\xi_\ell$ and detection efficiency $\eta = 0.9$. We calculated the estimator $\widehat{W}_{h,n}^\eta$ with optimal bandwidth $h_{\mathrm{opt}} = (\log n/(2\beta + 2\gamma))^{-1/2}$. We then reconsidered the kernel function and localized it by using a modified kernel having Fourier transform

$$\widetilde{K}_h^{\prime \eta}(t) = \frac{|t|}{2\tilde{N}^\eta(t/\sqrt{\eta})}$$

(19)

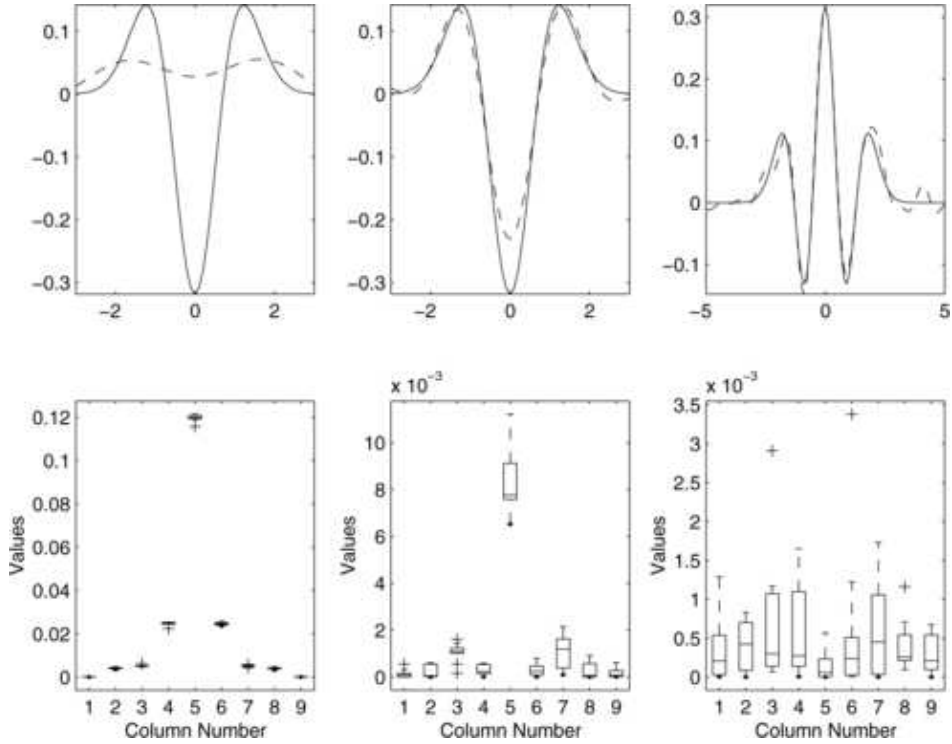

FIG. 2. *Left: One-photon state, $\eta = 0.9$, $n = 5000$. Middle: Same data, modified kernel. Right: Two-photon state, $\eta = 0.95$, $n = 10{,}000$.*

$$\times \left( I\left(|t| \leq \frac{1}{h}\right) + \exp\left(h^2 - \frac{1}{u(2/h - u)}\right) I\left(\frac{1}{h} \leq |t| \leq \frac{2}{h}\right) \right).$$

The function $\widetilde{K}'^\eta$ is an infinitely differentiable function (much smoother than $\widetilde{K}^\eta$); thus $K'^\eta$ decays exponentially fast. In Figure 2 we plot a transversal cut corresponding to the line $p = 0$, passing through the most difficult point to estimate $(0, 0)$, in which the error is dominated by the bias. The true Wigner function is plotted with a continuous line and the dashed line represents an estimator for one sample of size $n$. The graphics on the left-hand side concern the one-photon state with the original kernel estimator while the graphics in the middle show the estimator with the modified kernel (19) at the same bandwidth. An important improvement can be noticed in the case of the kernel $K'^\eta$. The left column concerns the two-photon state with modified kernel. The pointwise loss was then computed for ten samples (each of size $n = 5000$) at points $(0, 0), (0, \pm 0.5), (0, \pm 1), (0, \pm 1.5)$ and $(0, \pm 2)$ and the corresponding boxplots are shown in the lower panels of Figure 2. We notice that the highest losses are indeed observed at $(0, 0)$ and that the losses are quite stable from one sample to another. In the case of the Wigner function of the Schrödinger cat state we considered samples larger than 10,000 data which we binned in a $100 \times 100$ histogram. Figure 3 shows a contour plot of our estimator $\widehat{W}_{h,n}^{\eta}$ for a sample of size $n = 500,000$ and $\eta = 0.95$. Characteristic features are clearly visible: two Gaussian-shaped domes on the sides with positive (thin lines) and negative (thick lines) oscillations in the center. A similar estimator has been computed for $\eta = 0.85$ and Figure 4 shows different cuts through these estimators (dashed lines) compared with the true Wigner function (continuous line). The relatively worse performance in the case $\eta = 0.85$ is confirmed by Table 1 which gives the mean square errors over 100 samples of size $n$ at different peaked or flat points $(q, p)$ of the Wigner function and for the two different noise levels, $\eta = 0.95$ and $\eta = 0.85$. Tomographic reconstruction with real data was considered in [5]. However, in this reference no Gaussian deconvolution is performed. Thus one actually estimates a convoluted Wigner function $W_\rho^\eta = \mathcal{R}^{-1}[p_\rho^\eta]$ with usual parametric rate within $\log n$ factors [8]. We have tested such an estimator for the case of the Schrödinger cat state with $n = 500,000$ and $\eta = 0.85$ and cuts through the Wigner function. A result we obtained is the dash-dotted line shown in the panels in the lower part of Figure 4. This cut can be compared with the dashed line representing our estimator, which performs both inverse Radon transformation and deconvolution.

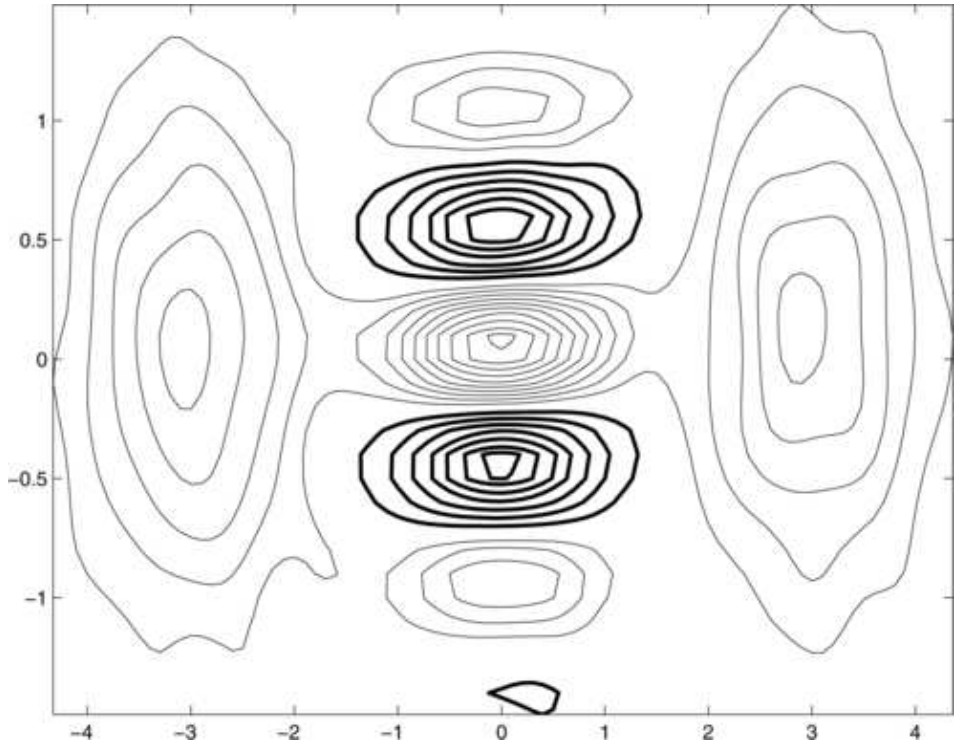

Fig. 3. *Contour plot of the estimated Wigner function for the Schrödinger cat state.*

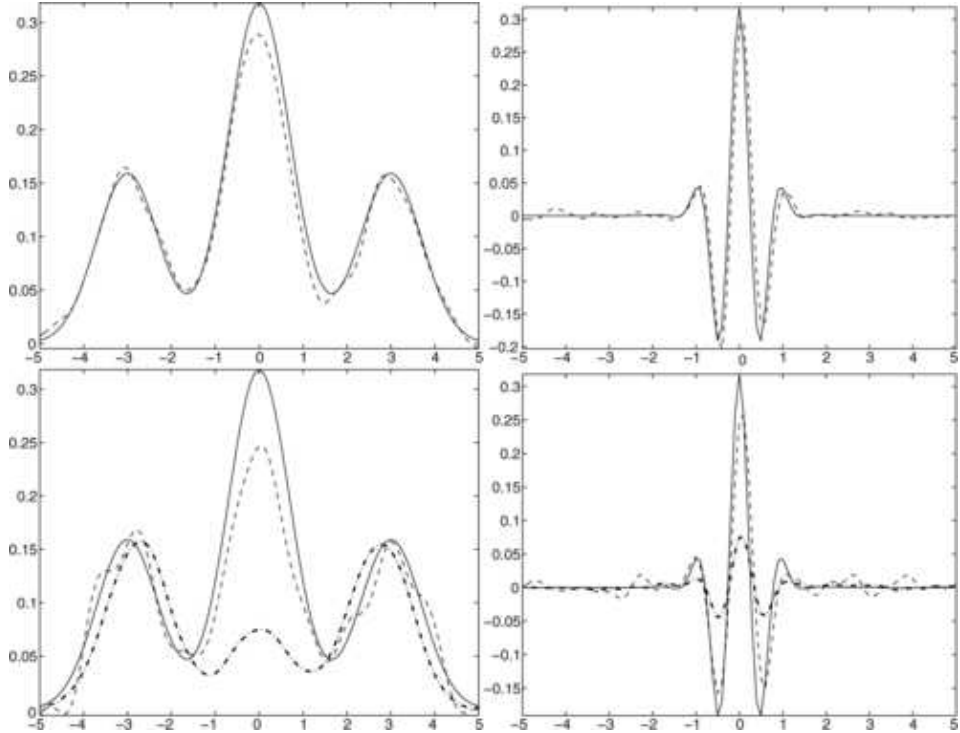

Fig. 4. *Transversal cuts through the Wigner function for the Schrödinger cat state.* Top: *Estimated Wigner function (dashed line),* $\eta = 0.95$, $n = 500,000$. Bottom: *Estimated Wigner function (dashed line) and estimator without deconvolution (dash-dotted line),* $\eta = 0.85$, *at* $n = 500,000$.

## 5. Proofs of upper bounds.

PROOF OF PROPOSITION 1. Since our data are i.i.d., we write

$$\mathbb{E}[\widehat{W}_{h,n}^{\eta}(z)] = \frac{1}{\pi}\int_0^{\pi}\int K_h^{\eta}([z,\phi] - y/\sqrt{\eta})p_{\rho}^{\eta}(y,\phi)\,dy\,d\phi$$

$$= \frac{1}{\pi}\int_0^{\pi} K_h^{\eta} * (\sqrt{\eta}p_{\rho}^{\eta}(\cdot\sqrt{\eta},\phi))([z,\phi])\,d\phi.$$

Now, write the convolution in the integral as an inverse Fourier transform. Indeed, it has Fourier transform [see (17)]

$$\mathcal{F}[K_h^{\eta} * (\sqrt{\eta}p_{\rho}^{\eta}(\cdot\sqrt{\eta},\phi))](t) = \widetilde{K}_h^{\eta}(t)\mathcal{F}_1[p_{\rho}^{\eta}(\cdot,\phi)](t/\sqrt{\eta})$$

$$= \tfrac{1}{2}|t|\mathcal{F}_1[p_{\rho}(\cdot,\phi)](t)I(|t| \le 1/h).$$

Replace this into the expected value of our estimator and use (7):

(20)
$$\mathbb{E}[\widehat{W}_{h,n}^{\eta}(z)] = \frac{1}{4\pi^2}\int_0^{\pi}\int_{-1/h}^{1/h} e^{-it[z,\phi]}|t|\widetilde{W}_{\rho}(t\cos\phi, t\sin\phi)\,dt\,d\phi$$

$$= \frac{1}{4\pi^2}\int\int e^{-i(qu+pv)}\widetilde{W}_{\rho}(u,v)I(\sqrt{u^2+v^2} \le 1/h)\,du\,dv$$

$$= \frac{1}{4\pi^2}\int e^{-i\langle z,w\rangle}\widetilde{W}_{\rho}(w)I(\|w\| \le 1/h)\,dw,$$

where we denote $w = (u,v)$. We recall that we also have

$$W_{\rho}(z) = \frac{1}{4\pi^2}\int e^{-i\langle z,w\rangle}\widetilde{W}_{\rho}(w)\,dw,$$

TABLE 1
*Schrödinger cat state: MSE $\times 10^5$ for 100 samples of size $n$ at points $(q,p)$ for $\eta = 0.95$ (left side) and $\eta = 0.85$ (right side)*

| $(q,p):n$ | 10,000 | 100,000 | 500,000 | 10,000 | 100,000 | 500,000 |
|---|---|---|---|---|---|---|
| $(0,0)$ | 507 | 173 | 119 | 1224 | 330 | 229 |
| $(0,3)$ | 54 | 10 | 4.16 | 428 | 161 | 67.9 |
| $(0,2.5)$ | 56.9 | 14.1 | 4.5 | 361 | 181 | 67.7 |
| $(0.5,0)$ | 414 | 113 | 70.1 | 909 | 258 | 164 |
| $(3,0)$ | 29.7 | 7.09 | 1.66 | 225 | 94.6 | 31.1 |

and then we write for the pointwise bias of our estimator,

$$|\mathbb{E}[\widehat{W}_{h,n}^\eta](z) - W_\rho(z)|^2 = \frac{1}{(4\pi^2)^2}\left|\int e^{-i\langle z,w\rangle}\{\mathcal{F}[E[\widehat{W}_{h,n}^\eta]](w) - \widetilde{W}_\rho(w)\}\,dw\right|^2$$

$$\leq \frac{1}{(4\pi^2)^2}\int |\widetilde{W}_\rho(w)|^2 e^{2\beta\|w\|^r}\,dw \int_{\|w\|>1/h} e^{-2\beta\|w\|^r}\,dw$$

$$\leq \frac{Lh^{r-2}}{4\pi\beta r}e^{-2\beta/h^r}(1+o(1)) \qquad \text{as } h\to 0,$$

by the assumption on our class. As for the variance of our estimator,

$$V[\widehat{W}_{h,n}^\eta(z)] = \mathbb{E}[|\widehat{W}_{h,n}^\eta(z) - \mathbb{E}[\widehat{W}_{h,n}^\eta(z)]|^2]$$

(21)
$$\leq \frac{1}{n}\mathbb{E}\left[\left|K_h^\eta\left([z,\Phi] - \frac{Y}{\sqrt{\eta}}\right)\right|^2\right]$$

$$\leq \frac{1}{n}\int_0^\pi \int (K_h^\eta([z,\phi] - y/\sqrt{\eta}))^2 p_\rho^\eta(y,\phi)\,dy\,d\phi.$$

At this point, let us denote

$$G(t) := \mathcal{F}[K_h^\eta([z,\phi] - \cdot/\sqrt{\eta})](t) = \sqrt{\eta}e^{it[z,\phi]\sqrt{\eta}}\widetilde{K}_h^\eta(-t\sqrt{\eta}).$$

Replace in (21) by taking into account that for a probability density $p_\rho^\eta(\cdot,\phi)$ we have $|\mathcal{F}_1[p_\rho^\eta(\cdot,\phi)]| \leq 1$,

$$\mathbb{E}\left[\left|K_h^\eta\left([z,\Phi] - \frac{Y}{\sqrt{\eta}}\right)\right|^2\right]$$

$$= \int_0^\pi \frac{1}{2\pi}\left|\int G*G(t)\mathcal{F}_1[p_\rho^\eta(\cdot,\phi)](t)\,dt\right|d\phi$$

$$\leq \frac{1}{2}\left(\int |G(t)|\,dt\right)^2 \leq \frac{1}{2}\left(\frac{\eta}{2}\int_{|t|\leq 1/(h\sqrt{\eta})}\frac{|t|}{\widetilde{N}^\eta(t)}\,dt\right)^2.$$

Finally we obtain

(22) $$\mathbb{E}\left[\left|K_h^\eta\left([z,\Phi] - \frac{Y}{\sqrt{\eta}}\right)\right|^2\right] \leq \frac{1}{2}\left(2\eta\int_0^{1/(h\sqrt{\eta})}\frac{t}{2}\exp\left(t^2\frac{1-\eta}{4}\right)dt\right)^2.$$

Let us note here that, more generally, for any positive $a, s$ and for any $A \in \mathbb{R}$, we can use integration by parts to get the asymptotic evaluation

(23) $$\int_0^x t^A \exp(at^s)\,dt = \frac{1}{as}x^{A+1-s}\exp(ax^s)(1+o(1)) \qquad \text{as } x\to\infty.$$

We use formula (23) for the integral in (22) as $1/h \to \infty$, and with (21) we get

$$V[\widehat{W}_{h,n}^\eta(z)] \leq \frac{2\eta^2}{(1-\eta)^2 n}\exp\left(\frac{1-\eta}{2\eta}\frac{1}{h^2}\right)(1+o(1)), \qquad n\to\infty.$$

□

PROOF OF THEOREM 2.   Over $\mathcal{B}$ we have

$$\mathbb{E}[|\widehat{W}_{h_{\mathrm{ad}},n}^{\eta}(z) - W_{\rho}(z)|^2] \leq \frac{L}{4\pi\beta r}(h_{\mathrm{ad}})^{r-2}\exp\left(-\frac{2\beta}{(h_{\mathrm{ad}})^r}\right)$$
$$+ \frac{2\eta^2}{(1-\eta^2)n}\exp\left(\frac{1-\eta}{2\eta(h_{\mathrm{ad}})^2}\right),$$

and it is easy to check that, for $(\beta, r, L) \in \mathcal{B}$,

$$\exp\left(-\frac{2\beta}{(h_{\mathrm{ad}})^r}\right) \leq \exp\left(-\frac{2\beta}{h_{\mathrm{opt}}^r}\right)(1 + o(1)),$$

$$\frac{1}{n}\exp\left(\frac{1-\eta}{2\eta(h_{\mathrm{ad}})^2}\right) = \exp\left(-\sqrt{\frac{\eta-1}{2\eta}\log n}\right) = o(1)\exp\left(-\frac{2\beta}{h_{\mathrm{opt}}^r}\right).$$

Thus, $\widehat{W}_{h_{\mathrm{ad}},n}^{\eta}$ attains precisely the rate $\varphi_n^2$ $(C = 1)$.   □

## 6. Proof of lower bounds.
In this section we will construct a pair of Wigner functions $W_1$ and $W_2$ depending on a parameter $\tilde{h}$ such that $\tilde{h} \to 0$ as $n \to \infty$. The choice of $\tilde{h}$ [see (31)] is such that it insures the existence of the lower bound in Theorem 1, and it should not be confused with the window $h$ appearing in the expression of the estimator which is optimal with respect to the upper bounds. We choose $W_1$ and $W_2$ of the forms

$$W_1(z) = W_0(z) + V_{\tilde{h}}(z) \quad \text{and} \quad W_2(z) = W_0(z) - V_{\tilde{h}}(z),$$

where $W_0$ is a fixed Wigner function corresponding to the density matrix $\rho_0$. The function $V_{\tilde{h}}$ is not a Wigner function of a density matrix but belongs to the linear span of the space of Wigner functions and thus has a corresponding matrix $\tau^{\tilde{h}}$ in the linear span of density matrices. The choice of $W_0, V_{\tilde{h}}$ is such that

$$\rho_1 = \rho_0 + \tau^{\tilde{h}} \quad \text{and} \quad \rho_2 = \rho_0 - \tau^{\tilde{h}}$$

are density matrices (positive and trace equal to 1) with Radon transforms $p_1$ and $p_2$. Suppose that the following conditions are satisfied:

(24)           $W_1$ and $W_2$       belong to the class $\mathcal{A}(\beta, r, L)$,

(25)           $|W_2(z) - W_1(z)| \geq 2\varphi_n(1 + o(1))$      as $n \to \infty$,

(26)
$$n\chi^2 := n\int_0^{\pi}\int\frac{(p_2^{\eta}(y,\phi) - p_1^{\eta}(y,\phi))^2}{p_1^{\eta}(y,\phi)}\,dy\,d\phi = o(1)$$

as $n \to \infty$.

Then we reduce the minimax risk to these two functions, $W_1$ and $W_2$, and bound the max from below by the mean of the two risks, to get for some $0 < \tau < 1$,

$$\inf_{\widehat{W}_n} \sup_{W_\rho \in \mathcal{A}(\beta, r, L)} \mathbb{E}[|\widehat{W}_n(z) - W_\rho(z)|^2]$$

$$\geq \left( \inf_{\widehat{W}_n} \frac{1}{2} \left( \mathbb{E}_{\rho_1}[|\widehat{W}_n(z) - W_1(z)|] \right. \right.$$

$$\left. \left. + (1 - \tau) \mathbb{E}_{\rho_1} \left[ I \left[ \frac{d\mathbb{P}_{\rho_2}^\eta}{d\mathbb{P}_{\rho_1}^\eta} \geq 1 - \tau \right] |\widehat{W}_n(z) - W_2(z)| \right] \right) \right)^2$$

$$\geq \frac{(1 - \tau)^2}{4} \cdot (2\varphi_n)^2 \mathbb{P}_{\rho_1}^2 \left[ \frac{d\mathbb{P}_{\rho_2}^\eta}{d\mathbb{P}_{\rho_1}^\eta} \geq 1 - \tau \right] (1 + o(1)).$$

We use the triangle inequality to get rid of the estimator and (25). Following Lemma 4 in [6], we know that the last probability in the display above is bounded from below by $1 - \tau^2$ provided that $n\chi^2 \leq \tau^4$. It is therefore sufficient to check (26), in order to find $\tau_n \to 0$, as $n \to \infty$ and give a lower bound of the minimax risk of order $\varphi_n^2(1 + o(1))$, for any estimator $\widehat{W}_n$.

We construct first the functions $W_{1,2}$ and then prove (24)–(26) in Section 6.3. Note that for the case $r = 2$ we prove a weaker form of (26): $n\chi^2 = O(1)$ as $n \to \infty$. The same reasoning as above shows that $\phi_n^2$ is then the optimal rate up to some constant (depending on some fixed $\tau$).

6.1. *Construction of the density matrix $\rho_0$.* In this section we will construct a family of density matrices $\rho^{\alpha, \xi}$ from which we will later select $\rho_0 = \rho^{\alpha_0, \xi_0}$ used in the lower bound. We derive their asymptotic behavior in Lemmas 1 and 2, and we show that $W_\alpha^\xi$ belongs to the class $\mathcal{A}(\beta, r, L)$ for $\alpha > 0$ small enough and $\xi$ close to 1.

Let us consider the Mehler formula (see [12], 10.13.22)

$$(27) \qquad \sum_{k=0}^\infty z^k \frac{1}{\sqrt{\pi} k! 2^k} H_k(x)^2 e^{-x^2} = \frac{1}{\sqrt{\pi(1 - z^2)}} \exp\left( -x^2 \frac{1 - z}{1 + z} \right),$$

where $H_k$ are the Hermite polynomials. Integrating both terms with $f_\alpha^\xi(z) = \alpha((1 - z)/(1 - \xi))^\alpha I(\xi \leq z \leq 1)$, for some $0 < \alpha, \xi < 1$, we get

$$(28) \qquad \begin{aligned} p_\alpha^\xi(x, \phi) &:= \sum_{k=0}^\infty \psi_k(x)^2 \int_0^1 f_\alpha^\xi(z) z^k \, dz \\ &= \int_0^1 \frac{f_\alpha^\xi(z)}{\sqrt{\pi(1 - z^2)}} \exp\left( -x^2 \frac{1 - z}{1 + z} \right) dz, \end{aligned}$$

where $\psi_k$ are the orthonormal vectors defined in (10). The Fourier transform of $p_\alpha^\xi$ is

$$(29) \qquad \widetilde{W}_\alpha^\xi(w) = \mathcal{F}_1[p_\alpha^\xi](\|w\|, \phi) = \int_0^1 \frac{f_\alpha^\xi(z)}{1-z} \exp\left(-\|w\|^2 \frac{1+z}{4(1-z)}\right) dz.$$

Notice that the normalization condition $\int p_\alpha^\xi = 1$ is equivalent to $\widetilde{W}_\alpha^\xi(0) = 1$, which is satisfied for the chosen functions $f_\alpha^\xi$, and thus $p_\alpha^\xi$ is a probability density. From the first equality in (28) we deduce that $p_\alpha^\xi$ is the probability density corresponding to a diagonal density matrix $\rho^{\alpha,\xi}$ with elements $\rho_{k,k}^{\alpha,\xi} = \int_0^1 z^k f_\alpha^\xi(z) \, dz$. We look now at the behavior of $p_\alpha^\xi(x, \phi)$ with respect to $x$.

LEMMA 1. *For all $0 < \alpha, \xi < 1$ and $|x| > 1$ there exist constants $c, C$ depending on $\alpha$ and $\xi$, such that $c|x|^{-(1+2\alpha)} \le p_\alpha^\xi(x, \phi) \le C|x|^{-(1+2\alpha)}$.*

PROOF. We have

$$p_\alpha^\xi(x, \phi) = \frac{\alpha}{(1-\xi)^\alpha \sqrt{\pi}} \int_\xi^1 \frac{(1-z)^{\alpha-1/2}}{(1+z)^{1/2}} \exp\left(-x^2 \frac{1-z}{1+z}\right) dz,$$

which by the change of variables $u = x\sqrt{\frac{1-z}{1+z}}$ becomes

$$p_\alpha^\xi(x, \phi) = \frac{\alpha 2^{\alpha+1} |x|}{(1-\xi)^\alpha \sqrt{\pi}} \int_0^{x\sqrt{(1-\xi)/(1+\xi)}} \frac{u^{2\alpha}}{(u^2 + x^2)^{\alpha+1}} \exp(-u^2) \, du.$$

By denoting $g(u) = u^{2\alpha} \exp(-u^2)$, the last integral is bounded for $|x| \ge 1$ as follows:

$$\frac{\alpha}{(1-\xi)^\alpha \sqrt{\pi} |x|^{2\alpha+1}} \int_0^{\sqrt{(1-\xi)/(1+\xi)}} g(u) \, du$$

$$\le p_\alpha^\xi(x, \phi) \le \frac{\alpha 2^{\alpha+1}}{(1-\xi)^\alpha \sqrt{\pi} |x|^{2\alpha+1}} \int_0^\infty g(u) \, du. \qquad \square$$

A similar analysis can be done for the matrix elements of $\rho^\alpha$. In the particular case $\alpha = 1$ and $\xi = 0$ we have $\rho_{k,k}^{1,0} = \frac{1}{(k+1)(k+2)}$.

LEMMA 2. *For all $0 < \alpha, \xi < 1$ we have*

$$\rho_{k,k}^{\alpha,\xi} = \frac{\alpha}{(1-\xi)^\alpha} \Gamma(\alpha+1) k^{-(1+\alpha)} (1 + o(1)) \qquad \textit{as } n \to \infty.$$

PROOF. We notice that by definition of $\rho_{k,k}^{\alpha,\xi}$ and the property

$$\int_0^1 z^k (1-z)^\alpha \, dz = \frac{\Gamma(1+\alpha)\Gamma(1+k)}{\Gamma(2+\alpha+k)},$$

$$\left| \rho_{k,k}^{\alpha,\xi} - \frac{\alpha}{(1-\xi)^\alpha} \frac{\Gamma(1+\alpha)\Gamma(1+k)}{\Gamma(2+\alpha+k)} \right| = \frac{\alpha}{(1-\xi)^\alpha} \int_0^\xi z^k (1-z)^\alpha \, dz$$

$$\leq \frac{\alpha \xi^{k+1}}{(1-\xi)^\alpha}.$$

Now, using Stirling's approximation for the function $\Gamma$ (see [1], 6.1.47) we deduce that $\Gamma(1+k)/\Gamma(2+\alpha+k) = k^{-(1+\alpha)}(1+o(1))$ and, given that $k^{1+\alpha}\xi^{k+1} = o(1)$, we obtain the desired result. $\quad \square$

LEMMA 3.  *For any* $(\beta, r, L)$ *such that* $0 < r \leq 2$, *there exist* $0 < \alpha, \xi \leq 1$ *such that* $W_\alpha^\xi$ *belongs to the class* $\mathcal{A}(\beta, r, L)$.

PROOF.  Using (29) we get

$$\int e^{2\beta \|w\|^r} |\widetilde{W}_\alpha^\xi(w)|^2 \, dw$$

$$= \int_0^\infty t e^{2\beta t^r} \left( \int_0^1 \frac{f_\alpha^\xi(z)}{1-z} \exp\left( -t^2 \frac{1+z}{4(1-z)} \right) dz \right)^2 dt$$

$$= \frac{\alpha^2}{(1-\xi)^{2\alpha}} \int_0^\infty t e^{2\beta t^r} \left( \int_\xi^1 (1-z)^{\alpha-1} \exp\left( -\frac{t^2}{2(1-z)} + \frac{t^2}{4} \right) dz \right)^2 dt$$

$$\leq \frac{\alpha^2}{(1-\xi)^{2\alpha}} \int_0^\infty t e^{2\beta t^r + t^2/2} \left( \int_\xi^1 (1-z)^{\alpha-1} \exp\left( -\frac{t^2}{2(1-\xi)} \right) dz \right)^2 dt$$

$$\leq \int_0^\infty t \exp\left( 2\beta t^r - \frac{t^2(1+\xi)}{2(1-\xi)} \right) dt \leq C(\beta, r, \xi),$$

where $C(\beta, r, \xi) > 0$ can be made smaller than $(2\pi)^2 L$ for any $0 < r \leq 2$ and for $0 < \xi < 1$ close enough to 1. $\quad \square$

6.2. *Construction of* $V_{\tilde{h}}$ *and asymptotic properties of* $\rho^{\tilde{h}}$. Let $V_{\tilde{h}}$ be the function defined on $\mathbb{R}^2$ whose Fourier transform is

(30)
$$\mathcal{F}_2[V_{\tilde{h}}](w) = \tilde{V}_{\tilde{h}}(w) := J_{\tilde{h}}(t)$$
$$= 2\sqrt{\pi \beta r L} \tilde{h}^{1-r/2} e^{\beta/\tilde{h}^r} e^{-2\beta |t|^r} J\left( |t|^r - \frac{1}{\tilde{h}^r} \right),$$

where $t = \|w\|$, and $J$ is a three-times continuously differentiable function with bounded derivatives and such that $I_{[2\delta, D-2\delta]}(u) \leq J(u) \leq I_{[\delta, D-\delta]}(u)$, for some $\delta > 0$ and $D > 4\delta$. The choice of the function $V_{\tilde{h}}$ is motivated for the case $0 < r < 2$ by the results on lower bounds for deconvolution obtained in [6]. The parameter $\tilde{h} \to 0$ as $n \to \infty$ is solution of the equation

(31)
$$\frac{2\beta}{\tilde{h}^r} + \frac{2\gamma}{\tilde{h}^2} = \log n + (\log \log n)^2.$$

When $r = 2$, we choose

$$(32) \qquad \tilde{h} = \left( \frac{\log(n \log n)}{2(\beta + \gamma)} \right)^{-1/2}.$$

We think of $V_{\tilde{h}}$ as a function belonging to the linear span of the Wigner functions. Indeed, as shown in (13), the convex map sending a density matrix $\rho$ to its corresponding Wigner function $W_\rho$ can be extended by linearity to an isometry (up to a constant) with respect to the $\|\cdot\|_2$ norm on the two spaces. We can thus construct a matrix $\tau^{\tilde{h}}$ belonging to the linear span of the space of density matrices and whose corresponding Wigner function is $V_{\tilde{h}}$. Because the function $V_{\tilde{h}}$ is invariant under rotations in the plane, the corresponding matrix has all off-diagonal elements equal to 0 and for the diagonal ones we can use the formula (from [21])

$$(33) \qquad \tau_{kk}^{\tilde{h}} = 4\pi^2 \int_0^\infty L_k(t^2/2) e^{-t^2/4} t J_{\tilde{h}}(t) \, dt,$$

where $L_k$ are the Laguerre polynomials defined in the proof of the following lemma.

LEMMA 4. *The matrix $\tau^{\tilde{h}}$ has the asymptotic behavior*

$$(34) \qquad \tau_{kk}^{\tilde{h}} = O(k^{-5/4}) o_{\tilde{h}}(1).$$

PROOF. We use the differential equation of the Laguerre polynomials (see [15], 8.979), $L_k(x) = \frac{1}{k}((x-1)L_k'(x) - x L_k''(x))$. Thus

$$\frac{d}{dt} L_k(t^2/2) = t L_k'(t^2/2) \quad \text{and} \quad \frac{d^2}{dt^2} L_k(t^2/2) = L_k'(t^2/2) + t^2 L_k''(t^2/2),$$

which implies

$$\frac{t^2}{2} L_k''(t^2/2) = \frac{1}{2} \frac{d^2}{dt^2} L_k(t^2/2) - \frac{1}{2} t^{-1} \frac{d}{dt} L_k(t^2/2)$$

and

$$L_k(t^2/2) = \frac{1}{2k} \left( (t^2 - 1) t^{-1} \frac{d}{dt} L_k(t^2/2) - \frac{d^2}{dt^2} L_k(t^2/2) \right).$$

Using integration by parts we obtain

$$\tau_{kk}^{\tilde{h}} = \frac{1}{k} \int_0^\infty L_k(t^2/2) e^{-t^2/4} [P_1(t) J_{\tilde{h}}(t) + P_2(t) J_{\tilde{h}}'(t) + P_3(t) J_{\tilde{h}}''(t)] \, dt,$$

with $P_i(t)$ polynomials of degree at most 3, whose coefficients do not depend on $\tilde{h}$ or $k$. As the support of the function under the integral is contained in the interval $[1/\tilde{h}, \infty)$, we can use the following bound for the behavior of Laguerre polynomials (see [27], Theorem 8.9.12): $\sup_{x \in [1,\infty)} e^{-x/2}|L_k(x)| = O(k^{-1/4})$. The matrix $\tau^{\tilde{h}}$ has thus the asymptotic behavior

$$\tau^{\tilde{h}}_{kk} \le Ck^{-5/4} \int_{1/\tilde{h}}^\infty |P_1(t)J_{\tilde{h}}(t) + P_2(t)J'_{\tilde{h}}(t) + P_3(t)J''_{\tilde{h}}(t)|\,dt$$

$$= O(k^{-5/4})o_{\tilde{h}}(1). \qquad \square$$

6.3. *Proofs of* (24)–(26) *involved in the lower bound.* Lemma 3 implies that for $\xi$ sufficiently close to 1, the Wigner function $W_\alpha^\xi$ belongs to the class $\mathcal{A}(\beta, r, a^2 L)$. On the other hand, combining the results of Lemma 2 and Lemma 4 we get that for any $\alpha < 1/4$ the diagonal matrices $\rho_1 = \rho^{\alpha,\xi} + \tau^{\tilde{h}}$ and $\rho_2 = \rho^{\alpha,\xi} - \tau^{\tilde{h}}$ are positive and have trace 1 for $\tilde{h}$ sufficiently small. Thus there exist $\alpha_0, \xi_0$ such that the corresponding $\rho_1$ and $\rho_2$ are density matrices and $W_0 = W_{\alpha_0}^{\xi_0} \in \mathcal{A}(\beta, r, a^2 L)$.

In the following proofs the constants $\delta$ and $D$ appear from the construction of $V_{\tilde{h}}$. The whole proof holds for arbitrarily small $\delta > 0$ and arbitrarily large $D > 4\delta$, hence the desired results.

PROOF OF (24). By the triangle inequality

$$\|\mathcal{F}_2[W_{1,2}]e^{\beta\|\cdot\|^r}\|_2 \le \|\mathcal{F}_2[W_0]e^{\beta\|\cdot\|^r}\|_2 + \|\mathcal{F}_2[V_{\tilde{h}}]e^{\beta\|\cdot\|^r}\|_2.$$

The first term in the sum above is less than $2\pi\sqrt{L}a$. For the second one we have

$$\int |\mathcal{F}_2[V_{\tilde{h}}](w)|^2 e^{2\beta\|w\|^r}\,dw$$

$$= \int_0^\pi \int |t||\mathcal{F}_2[V_{\tilde{h}}](t\cos\phi, t\sin\phi)|^2 e^{2\beta|t|^r}\,dt\,d\phi$$

$$= \pi \int |t||J_{\tilde{h}}(t)|^2 e^{2\beta|t|^r}\,dt$$

$$\le 4\pi^2 \beta r L \tilde{h}^{2-r} e^{2\beta/\tilde{h}^r} \int_{\delta \le |t|^r - 1/\tilde{h}^r \le D - \delta} |t| e^{-2\beta|t|^r}\,dt$$

$$\le 4\pi^2 L e^{-2\beta\delta}.$$

Thus, if we take $a = 1 - e^{-\beta\delta/2}$, we get $W_{1,2}$ in the class $\mathcal{A}(\beta, r, L(1 - e^{-\beta\delta/2} + e^{-\beta\delta}))$ included in $\mathcal{A}(\beta, r, L)$. $\square$

PROOF OF (25). Notice that $|W_2(z) - W_1(z)|^2$ is equal to

$$\left|\frac{1}{4\pi^2} \int_{\mathbb{R}^2} e^{-i\langle z, w\rangle} (\widetilde{W}_2(w) - \widetilde{W}_1(w))\,dw\right|^2$$

$$= \left| \frac{1}{4\pi^2} \int_0^{2\pi} \int_0^{\infty} e^{-it[z,\phi]} |t| (\widetilde{W}_2(t\cos\phi, t\sin\phi) \right.$$
$$\left. - \widetilde{W}_1(t\cos\phi, t\sin\phi)) \, dt \, d\phi \right|^2$$
$$= \left| \frac{1}{2\pi^2} \int_0^{2\pi} \int_0^{\infty} e^{-it[z,\phi]} |t| J_{\bar{h}}(t) \, dt \, d\phi \right|^2.$$

Take $z = 0$ without loss of generality:

$$|W_2(z) - W_1(z)|^2$$
$$= \left| \frac{1}{2\pi} \int_0^{\pi} \int |t| J_{\bar{h}}(t) \, dt \right|^2$$
$$\geq 4\pi\beta r L \tilde{h}^{2-r} e^{2\beta/\tilde{h}^r} \left| \frac{1}{2\pi} \int_{2\delta \leq |t|^r - 1/\tilde{h}^r \leq D - 2\delta} |t| e^{-2\beta|t|^r} \, dt \right|^2$$
$$\geq 4 \frac{L}{4\pi\beta r} \tilde{h}^{r-2} e^{-2\beta/\tilde{h}^r} [e^{-4\beta\delta}(1 + o(1)) - e^{-2\beta(D-2\delta)}(1 + o(1))]^2,$$

which is larger than $4\varphi_n^2 [e^{-4\beta\delta} - e^{-2\beta(D-2\delta)}]^2 (1 + o(1))$ for $n$ large enough. Note that for $0 < r < 2$, the $\tilde{h}$ solution of (31) provides exact lower bounds, while for $r = 2$, $\bar{h}$ given by (32) provides optimal rates of order $(n \log n)^{-\beta/(\beta+\gamma)}$, which are within a logarithmic factor optimal. $\square$

PROOF OF (26). We want to bound from above $n\chi^2 \leq \pi n \int (p_2^{\eta}(y) - p_1^{\eta}(y))^2 / p_1^{\eta}(y) \, dy$. We have proven that $p_1(x) \geq C x^{-2}$ for all $|x| \geq 1$. It is easy to prove that after convolution with the Gaussian density of the noise the asymptotic decay cannot be faster; thus $p_1^{\eta}(y) \geq \frac{c_1}{y^2}, \forall |y| \geq M$, for some fixed $M > 0$. Then we split the integration domain into $|y| \leq M$ and $|y| > M$ and get

$$(35) \qquad n\chi^2 \leq Cn\left( C(M) \|p_2^{\eta} - p_1^{\eta}\|^2 + \int_{|y|>M} y^2 (p_2^{\eta}(y) - p_1^{\eta}(y))^2 \, dy \right).$$

Let us see first that

$$(36) \qquad \begin{aligned} \|p_2^{\eta} - p_1^{\eta}\|^2 &= C \int |J_{\tilde{h}}(t)|^2 e^{-(1-\eta)t^2/(2\eta)} \, dt \\ &\leq C\tilde{h}^{1-r} \exp\left(\frac{2\beta}{\tilde{h}^r}\right) \int_{(1+\delta\tilde{h}^r)^{1/r}/\tilde{h}}^{\infty} e^{-4\beta t^r - (1-\eta)t^2/(2\eta)} \, dt \\ &\leq C\tilde{h}^{2-r} \exp\left(-\frac{2\beta}{\tilde{h}^r} - \frac{1-\eta}{2\eta\tilde{h}^2}\right). \end{aligned}$$

Then

$$\int_{|y|>M} y^2 (p_2^{\eta}(y) - p_1^{\eta}(y))^2 \, dy$$

$$\text{(37)} \quad \leq \int \left( \frac{\partial}{\partial t} (J_{\tilde{h}}(t) e^{-(1-\eta)t^2/(4\eta)}) \right)^2 dt$$

$$\leq C \tilde{h}^{1-r} \exp\left( \frac{2\beta}{\tilde{h}^r} \right) \int_{(1+\delta \tilde{h}^r)^{1/r}/\tilde{h}}^{\infty} t^2 e^{-4\beta t^r - (1-\eta)t^2/(2\eta)} \, dt$$

$$\leq C \tilde{h}^{-r} \exp\left( -\frac{2\beta}{\tilde{h}^r} - \frac{1-\eta}{2\eta \tilde{h}^2} \right).$$

For the case $0 < r < 2$ choose $\tilde{h}$ as solution of (31) to get that the expressions in (36) and (37) tend to 0, and together with (35) this concludes the proof of (24). For the case $r = 2$, $\tilde{h}$ given by (32), we get that the expression in (36) tends to 0 and (37) stays bounded as $n \to \infty$; thus we obtain the desired result. □

C. Butucea                                M. Guţă
Modal'X                                   School of Mathematical Sciences
Université Paris X                        University of Nottingham
200, avenue de la République              Nottingham NG7 2RD
92001 Nanterre Cedex                      United Kingdom
France                                    E-mail: madalin.guta@nottingham.ac.uk
and
PMA
175, rue de Chevaleret
75013 Paris
France
E-mail: butucea@ccr.jussieu.fr

                    L. Artiles
                    Eurandom
                    P.O. Box 513
                    5600 MB Eindhoven
                    The Netherlands
                    E-mail: artiles@eurandom.tue.nl